\documentclass[12pt]{article}
\def\hind{\hangindent=2pc\hangafter=1}
\newfont{\smcaps}{cmcsc10 scaled\magstep1}

\newcommand{\Cov}{{\rm ~Cov\,}}

\newcommand{\AIC}{{\rm ~AIC\,}}
\newcommand{\BIC}{{\rm ~BIC\,}}

\newcommand{\AR}{{\rm ~AR\,}}
\newcommand{\ARphi}{{\rm ~AR\,}_\phi}
\newcommand{\ARzeta}{{\rm ~AR\,}_\zeta}
\newcommand{\AICzeta}{{\rm ~AIC\,}_\zeta}
\newcommand{\BICzeta}{{\rm ~BIC\,}_\zeta}

\newcommand{\EstSd}{{\rm ~EstSd\,}}
\font\tty=cmtt10 at 12truept
\usepackage{epsfig}
\raggedright
\begin{document}

\title{SUBSET AUTOREGRESSION: A NEW APPROACH}
\date{}
\author{{\smcaps By A.I. McLeod AND Y. Zhang}\\
{\it The University of Western Ontario and  Acadia University\/}
}
\maketitle
\bigskip
\hrule
\bigskip

\noindent A. Ian McLeod and Ying Zhang, (2006).  
Partial autocorrelation parameterization for subset autoregression.  
{\it Journal of Time Series Analysis} 27/4, 599-612.
\smallskip
\noindent 10.1111/j.1467-9892.2006.00481.x

\newpage
\baselineskip=22pt
\hrule
\bigskip
{\bf Abstract.\/}
A new version of the partial autocorrelation plot and
a new family of subset autoregressive models are introduced.
A comprehensive approach to model identification, estimation and diagnostic checking
is developed for these models.
These models are better suited to efficient model building of high-order autoregressions
with long time series.
Several illustrative examples are given.

\hfill\eject

\begin{center}
1. INTRODUCTION
\end{center}
\medskip

The AR$(p)$ model with mean $\mu$ may be written,
$\phi(B)(z_t-\mu) = a_t$, where $\phi(B)=1-\phi_1 B-\cdots \phi_p B^p$.
$B$ is the backshift operator on $t$ and $a_t, t=1,\ldots,n$ are normal and independently
distributed with mean zero and variance $\sigma_a^2$.
The admissible region for stationary-causal autoregressive processes is defined
by the region for which all roots of $\phi(B)=0$ lie outside the
unit circle (Brockwell and Davis, 1991).
The usual subset autoregressive model is obtained by constraining some of the
$\phi$-parameters to zero. In this case we may write,
$\phi(B)=1-\phi_{i_1} B^{i_1}-\cdots \phi_{i_m} B^{i_m}$
where $i_1<\ldots<i_m$.
This model will be denoted by $\ARphi(i_1,\ldots,i_m)$.
Such subset autoregressive models are often used for modelling seasonal or periodic
time series as well as for obtaining a more parsimonious representation
of autoregressive models.
McClave (1975) presented an algorithm using the Yule-Walker estimators which may
be used to find the best model according to some criterion such as the AIC or BIC.
However, as pointed out by Haggan and Oyetunji (1983), the algorithm given by McClave (1975)
only finds the single best solution although in practice it is often desirable to
examine a range of plausible models.
The algorithm developed given by Haggan and Oyetunji (1983), utilizing least squares,
is not as computationally efficient as that of McClave (1983) but it is easy find the
best set of models.
One drawback of this approach is that it is based on least squares.
Least squares estimates may be preferable to Yule-Walker estimates due to their lower
bias but least squares estimates may be inadmissible.
Although an admissible model may not be needed for short-term forecasting,
it is required for spectral estimation or data simulation in engineering design
(Hipel and McLeod, 1994, \S 9.7.3).
Zhang and Terrell (1997) have suggested a new criterion, the projection modulus,
which is computationally more efficient but their method is based on
Yule-Walker estimates that are known to be less accurate than some
alternatives (Tj\o stheim and Paulsen, 1983; Percival and Walden, A.T, 1993, p.414 and p.453;
Zhang and McLeod, 2005).
Bayesian methods of variable selection in regression were introduced by George
and McCulloch (1993) and Bayesian methods for subset autoregression have been
developed by Chen (1999) and Unnikrishnan (2004).
The approach developed in this paper is computationally more efficient than previous techniques,
based on maximum likelihood
and
well suited to fitting long time series and high dimensional subset
autoregressive models as is illustrated in \S 3.3.

We now introduce the new subset autoregression models.
Consider the Durbin-Levinson recursion
\begin{equation}
\phi_{j,k+1}=\phi_{j,k}-\phi_{k+1,k+1} \phi_{k+1-j,k},\quad j=1,\ldots,k
\newcounter{DLRecursion}
\setcounter{DLRecursion}{\value{equation}}
\end{equation}
where $k=1,\ldots,p$ and $\zeta_i=\phi_{i,i}$.
This recursion can be used to define a transformation,
\begin{equation}
{\cal B}: (\zeta_1,\ldots,\zeta_p)\rightarrow (\phi_1,\ldots,\phi_p),
\newcounter{BTransformation}
\setcounter{BTransformation}{\value{equation}}
\end{equation}
that is one-to-one, continuous and differentiable
inside the admissible region (Barndorff-Neilsen and Schou, 1973).
Both $\cal B$ and its inverse ${\cal B}^{-1}$ are easily computed (Monahan, 1984).
To extend this transformation to the subset autoregressive case we simply constrain
some of the $\zeta$-parameters to zero.
In general, this subset AR model may be denoted by $\ARzeta(i_1,\ldots,i_m)$ where
the underlying parameters are $(\zeta_{i_1},\ldots,\zeta_{i_m})$.
The $\ARphi(i_1,\ldots,i_m)$ and $\ARzeta(i_1,\ldots,i_m)$ are similar but distinct models.
For example, in the $\ARphi(1,3)$,
$\zeta_1=\phi_1/(1-\phi_1 \phi_3 - \phi_3^2), \zeta_2=\phi_1 \phi_3 / (1-\phi_3^2)$
and $\zeta_3=\phi_3$,
whereas in the $\ARzeta(1,3)$,
$\phi_1=\zeta_1$, $\phi_2=-\zeta_1 \zeta_3$ and $\phi_3=\zeta_3$.

In general, the $\ARzeta(i_1,\ldots,i_m)$ model
has only $m$ parameters, not including $\mu$ and $\sigma_a^2$,
but it specifies $p=i_m$ values for the parameters in $\phi$-space,
$\phi_1(\zeta_1,\ldots,\zeta_{i_m}),\ldots,\phi_p(\zeta_1,\ldots,\zeta_{i_m})$.
In $\phi$-space the admissible region is
a complex $m$-dimensional subspace of the original $p$-dimensional space
of $(\phi_1,\ldots,\phi_p)$.
In contrast, the admissible region in the $\zeta$-space,
${\cal D}_\zeta$,
for the $\ARzeta(i_1,\ldots,i_m)$ model is simply the $m$
dimensional cube with boundary surfaces corresponding to $\pm 1$.
The transformation ${\cal B}$ induces the transformation
${{\cal B}_{i_1,\ldots,i_m}}: (\zeta_{i_1},\ldots,\zeta_{i_m})\rightarrow (\phi_1,\ldots,\phi_p)$
defined by setting $\zeta_i = 0$ for $i \not \in {i_1,\ldots,i_m}$ in (\theBTransformation).
Denote the image of ${\cal D}_\zeta$ using the transformation
${{\cal B}_{i_1,\ldots,i_m}}$
by ${\cal D}_\phi$.
Then ${\cal D}_\phi$ is a very complicated subset of the original $p$ dimensional
admissible space of the full $\AR(p)$ model.
From Barndorff-Neilsen and Schou (1973, Theorem 2) it follows that the transformation
${\cal B}_{i_1,\ldots,i_m}: {\cal D}_\zeta \rightarrow {\cal D}_\phi$
is one-to-one, continuous and differentiable.
Denote the $p$ functions determined by ${{\cal B}_{i_1,\ldots,i_m}}$
as $\phi_j(\zeta_{i_1},\ldots,\zeta_{i_m}),\ j=1,\ldots,p$.
It follows from (\theDLRecursion) that each of these $p$ functions,
$\phi_1(\zeta_{i_1},\ldots,\zeta_{i_m}),\ldots,\phi_p(\zeta_{i_1},\ldots,\zeta_{i_m})$,
are polynomial functions of $\zeta_{i_1},\ldots,\zeta_{i_m}$.

\begin{center}
2. MODEL FITTING
\end{center}

{\noindent \it 2.1 Exact likelihood function\hfill}

The sample mean is asymptotically efficient so we will assume the series $z_1,\ldots,z_n$
after mean correction has mean zero.
Then the exact loglikelihood function, apart from a constant, for the
$\ARzeta(\zeta_{i_1},\ldots,\zeta_{i_m})$ may be written,
\begin{equation}
{\cal L}(\zeta, \sigma_a^2)=-{1\over 2}\log(\det(\Gamma_n))-{1\over 2}z^{\prime} \Gamma_n^{-1} z
\newcounter{Loglikelihood}
\setcounter{Loglikelihood}{\value{equation}}
\end{equation}
where $z=(z_1,\ldots,z_n)$ and $\Gamma_n$ is the $n \times n$ covariance matrix
with $(i,j)$ entry $\gamma_{i-j}=\Cov(z_{t-i},z_{t-j})$.
It follows from Box Jenkins and Reinsel (1994, \S A7.5),
$z^{\prime} \Gamma_n^{-1} z ={\cal S}(\zeta)/\sigma_a^2$,
where ${\cal S}(\zeta)=\beta^{\prime} D \beta$,
$\beta=(-1,\phi_1(\zeta),\ldots,\phi_p(\zeta))$
and
$D$ is the $(p+1) \times (p+1)$ matrix with $(i,j)$-entry,
$D_{i,j}=z_i z_j + \cdots + z_{n-j} z_{n-i}$.
Then letting $p=i_m$
and
$g_p=\det(\Gamma_n/\sigma_a^2)=\det(\Gamma_p/\sigma_a^2),$
we have from Barndorff-Neilsen and Schou (1973, eqns. 5 and 8),
\begin{equation}
g_p=\prod_{i\in \{i_1,\ldots,i_m \}} (1-\zeta_i^2)^{-i}.
\newcounter{CovarianceDeterminant}
\setcounter{CovarianceDeterminant}{\value{equation}}
\end{equation}
The loglikelihood function can now be written,
\begin{equation}
{\cal L}(\zeta,\sigma_a^2)=-{n\over 2}\log(\sigma_a^2)- {1\over 2}\log(g_p)
- {1 \over 2 \sigma_a^2} {\cal S}(\zeta).
\newcounter{FinalLoglikelihood}
\setcounter{FinalLoglikelihood}{\value{equation}}
\end{equation}
Maximizing over $\sigma_a^2$ and dropping constant terms, the concentrated loglikelihood is,
\begin{equation}
{\cal L}_c(\zeta)=-{n\over 2}\log(\hat \sigma_a^2)- {1\over 2}\log(g_p),
\newcounter{ConcentratedLoglikelihood}
\setcounter{ConcentratedLoglikelihood}{\value{equation}}
\end{equation}
where $\hat \sigma_a^2 = {\cal S}(\zeta)/n$.
${\cal L}_c$ can be optimized numerically using a constrained optimization algorithm
such as {\tty FindMinimum\/} in {\it Mathematica\/}.
The initial evaluation of $D$ requires $O(n)$ flops but this is only done once
and subsequent likelihood evaluations only require $O(p^2)$ flops.

There are a number of algorithms for ARMA likelihood evaluation
and many of these are listed in (Box and Luce\~no, 1997, \S 12B).
Anyone of these algorithms could also be used.
However, all of these algorithms require $O(n)$ flops per likelihood evaluation
whereas the algorithm given in this section only requires $O(p^2)$
and so is much more efficient.

As shown in Theorem 2 in \S 2.2,
statistically efficient initial values of the parameters may be obtained
using the partial autocorrelations computed by the Burg algorithm.
With this approach it is possible to obtain exact maximum likelihood
estimates for even quite high-order AR models as illustrated by
the monthly sunspot example, \S 3.3, where $m=70$ coefficients were estimated.

{\noindent \it 2.2 Large-sample distribution of the estimates\hfill}

For an observed time series $z_1,\ldots,z_n$ generated by an $\ARzeta(i_1,\ldots,i_m)$
model, let $\hat \zeta = (\hat \zeta_{i_1}, \ldots, \hat \zeta_{i_m})$ denote the
maximum likelihood estimate of $\zeta = (\zeta_{i_1}, \ldots, \zeta_{i_m})$.
In the full model case, $m=p$ and $(\zeta_1, \ldots, \zeta_p)$ is a reparameterization
of the $\AR(p)$ model.
However, in the subset case when $m<p$,
the parameters $\phi_1,\ldots,\phi_p$ are constrained
and so $\hat \phi_1,\ldots,\hat \phi_p$
do not have the usual distribution due to these constraints.
The following theorem is established in Appendix A.

{\smcaps Theorem 1.\/}
$\hat{\zeta} \stackrel{p}{\longrightarrow} \zeta$
and
$\surd n (\hat \zeta - \zeta)
\stackrel{D}{\longrightarrow} {\cal N}(0,{\cal I}_\zeta^{-1})$,
where
$\stackrel{p}{\longrightarrow}$ denotes convergence in probability,
$\stackrel{D}{\longrightarrow}$ denotes convergence in distribution and
${\cal I}_\zeta$ is the large-sample Fisher information matrix per observation of $\zeta$.

Properties of the information matrix ${\cal I}_\zeta$ are discussed in
Barndorff-Neilsen and Schou (1973)
for the case of the full model, $\ARzeta(1,\ldots,p)$, but a general method
of computing ${\cal I}_\zeta$ is not explicitly given.
It is shown in Appendix A that
${\cal I}_\zeta={\cal J}_\zeta^{\prime} {\cal I}_\phi {\cal J}_\zeta,$
where
\begin{equation}
{\cal J}_\zeta ={\partial (\phi_1,\ldots,\phi_p) \over \partial (\zeta_{i_1},\dots,\zeta_{i_m})}
\newcounter{JacobianZeta}
\setcounter{JacobianZeta}{\value{equation}}
\end{equation}
and ${\cal I}_\phi$ is the information matrix for
$\phi_1,\ldots,\phi_p$ in the unrestricted $\AR(p)$ model.
Since ${\cal I}_\phi=\sigma_a^{-2} \Gamma_p$, ${\cal I}_\phi$
may be easily computed using the result of Siddiqui (1958),
\begin{equation}
{\cal I}_\phi =
\left(
\sum\limits_{k=0}^{\min(i,j)} \phi_{i-k+1} \phi_{j-k+1}-\phi_{p+k-i+1} \phi_{p+k-j+1}
\right)_{p \times p},
\newcounter{SiddiquiFormula}
\setcounter{SiddiquiFormula}{\value{equation}}
\end{equation}
where $\phi_0 = -1$.
The Jacobian ${\cal J}_\zeta$ is quite complicated.
First, consider the full model case, $\ARzeta(1,\ldots,p)$.
The required Jacobian may be derived as the product of a sequence
of Jacobians of transformations used in the Durbin-Levinson
algorithm, eq. (\theDLRecursion), to obtain,
${\cal J}_\zeta={\cal J}_{p-1} \cdots {\cal J}_{1}$,
where
\begin{equation}
{\cal J}_{k}=
{\partial ( \phi_{k+1,1},\ldots,\phi_{k+1,p} ) \over
\partial ( \phi_{k,1},\ldots,\phi_{k,p} ) },
\newcounter{Jacobiank}
\setcounter{Jacobiank}{\value{equation}}
\end{equation}
where $\phi_{k,j}=\phi_{k,k}$ for $j\ge k$ and otherwise
for $j<k$, $\phi_{k,j}$ is as defined in (\theDLRecursion).
It may then be shown that
\begin{equation}
{\cal J}_{k}=\left (
\begin{array}{cc}
J_{p-k} & A_{p-k,k} \\
0_{k,p-k} & I_k
\end{array}
 \right ),
\newcounter{JacobiankMatrix}
\setcounter{JacobiankMatrix}{\value{equation}}
\end{equation}
where $J_{p-k}$ is the $(p-k) \times (p-k)$ matrix with $(i,j)$-entry,
$J_{i,j}$, where
\begin{equation}
J_{i,j}=\left\{
\begin{array}{cl}
1 & {\rm if}\quad i=j\\
-\zeta_{p-k+1} & {\rm if}\quad i=p-k+1-j \wedge i \ne j \\
1-\zeta_{p-k+1} & {\rm if}\quad i=p-k+1-j \wedge i = j \\
0 & {\rm otherwise}
\end{array}
\right. ,
\newcounter{Jij}
\setcounter{Jij}{\value{equation}}
\end{equation}
$A_{p-k,k}$ is the $(p-k) \times k$ matrix whose first column is
$(-\phi_{p-k,p-k},-\phi_{p-k,p-k-1},\ldots, -\phi_{p-k,1})$
and whose remaining elements are zeros,
$0_{p-k,k}$ is the $k \times (p-k)$ matrix with all
zero entries,
and $I_k$ is the $k \times k$ identity matrix.
For example, for the $\ARzeta(1,2,3,4)$,
\begin{eqnarray}
{\cal J}_\zeta & = &
\left(
\begin{array}{cccc}
1         & 0         & -\zeta_4 & -\zeta_3 \\
0         & 1-\zeta_4 & 0        & -\zeta_2 - \zeta_1(1+\zeta_2)\zeta_3 \\
-\zeta_4  & 0         & 1        & -\zeta_1(1+\zeta_2) - \zeta_2 \zeta_3 \\
0         & 0         & 0        & 1
\end{array}
\right)
\nonumber\\
&&
\times
\left(
\begin{array}{cccc}
1-\zeta_2 & -\zeta_1 & 0 & 0 \\
0         & 1        & 0 & 0 \\
0         & 0        & 1 & 0 \\
0         & 0        & 0 & 1
\end{array}
\right)
\times
\left(
\begin{array}{cccc}
1         & -\zeta_3 & -\zeta_2            & 0 \\
-\zeta_3  & 1        & -\zeta_1(1+\zeta_2) & 0 \\
0         & 0        & 1                   & 0 \\
0         & 0        & 0                   & 1
\end{array}
\right).
\newcounter{Jfour}
\setcounter{Jfour}{\value{equation}}
\end{eqnarray}

The information matrix of $\zeta$ in the subset case,
$\ARzeta(i_1,\ldots,i_m)$,
may be obtained by selecting rows and columns corresponding
to $i_1,\ldots,i_m$ from full model information matrix.
Equivalently, the information matrix in the subset case could
also be obtained by selecting the columns corresponding
to $i_1,\ldots,i_m$ in the Jacobian matrix corresponding to the full model case
to obtain ${\cal J}_{(\zeta_{i_1},\ldots,\zeta_{i_m})}$.
Then
\begin{equation}
{\cal I}_{(\zeta_{i_1},\ldots,\zeta_{i_m})} =
{\cal J}_{(\zeta_{i_1},\ldots,\zeta_{i_m})}^{\prime}
{\cal I}_\phi
{\cal J}_{(\zeta_{i_1},\ldots,\zeta_{i_m})}.
\newcounter{InformationMatrixSubsetZeta}
\setcounter{InformationMatrixSubsetZeta}{\value{equation}}
\end{equation}
For example, using our {\it Mathematica\/} software, for the $\ARzeta(1,12)$,
\begin{equation}
{\cal I}_{\zeta}=
\frac{1}{1 - \zeta_{12}^2}
\left(
\begin{array}{cc}
\left( 1 - 2 \zeta_1^{10} \zeta_{12} + \zeta_{12}^2 \right) /
\left( 1 - \zeta_1^2 \right)
& 0 \\ 0 & 1
\end{array}
\right).
\newcounter{InformationMatrixSeasonalExample}
\setcounter{InformationMatrixSeasonalExample}{\value{equation}}
\end{equation}

As a check on the formula for ${\cal I}_\zeta$,
an $\ARzeta(1,2,3,4)$ with $\zeta_i=0.5,\ i=1,\ldots,4$ and $n=500$
was simulated and fit 1,000 times.
The empirical covariance matrix of the $\zeta$-parameters was found to agree closely
with the theoretical covariance matrix $n^{-1} {\cal I}_\zeta^{-1}$.
This experiment was repeated using the $\ARzeta(1,4)$ model.

As mentioned in \S 2.1, the Burg algorithm can be used to generate good
initial estimates of the parameters $\zeta_i, i=i_1,\ldots,i_m$.
As shown in Theorem 2 below these estimates are asymptotically fully
efficient.
However, we prefer to use the exact MLE for our final model
estimates since these estimates are known to be second order efficient
(Taniguchi, 1983)
and simulation experiments have shown that the exact MLE
estimates usually perform better than alternatives in small samples
(Box and Luce\~no, 1997, \S 12B).

{\smcaps Theorem 2.\/} In the $\ARzeta(i_1,\ldots,i_m)$ model
let $\hat{\phi}_{k,k}$, $k=1,\ldots,p$, where $p \ge i_m$
denote the partial autocorrelations estimated using the Burg algorithm.
Then $\hat{\phi}_{i_1,i_1},\ldots,\hat{\phi}_{i_m,i_m}$
are asymptotically efficient estimates for $\zeta_{i_1},\ldots,\zeta_{i_m}$.

Theorem 2 follows from the fact that the Burg algorithm provides asymptotically efficient
estimates (Percival and Walden, 1993, p.433)
of $\zeta_1, \ldots, \zeta_p$ in the full $\AR(p)$ model.
Then the large-sample covariance matrix of
$\hat{\phi}_{i_1,i_1},\ldots,\hat{\phi}_{i_m,i_m}$,
given by eq. (\theInformationMatrixSubsetZeta),
is seen to be the same as that $\hat{\zeta}_{i_1},\ldots,\hat{\zeta}_{i_m}$.

{\noindent \bf 2.3 Model identification\hfill}

Theorem 1 provides the basis for a useful model identification method for $\ARzeta(i_1,\ldots,i_m)$
using the partial autocorrelation function.
The partial autocorrelations are estimated for a suitable number of lags $k=1,\ldots, K$.
Typically $K<n/4$.
We recommend the Burg algorithm for estimating the partial autocorrelations
$\hat \phi_{k,k}=\hat{\zeta}_{k}$
since it provides more accurate estimates of the partial autocorrelations in many situations
(Percival and Walden, 1993, p.414) than the Yule-Walker algorithm.
Based on the fitted $\AR(K)$ model the estimated standard errors, $\EstSd(\hat \zeta_k)$,
of $\hat \zeta_k$ are obtained.
A suitable $\ARzeta(i_1,\ldots,i_m)$ model may be selected by examining a plot
of
$\hat \zeta_k \pm 1.96 \EstSd(\hat \zeta_k)$ vs $k$.
This modified partial autocorrelation plot is generally more
useful than the customary one (Box, Jenkins and Reinsel, 1994).
The use of this partial autocorrelation plot is illustrated in \S 3.1$-$3.2.

Another method of model selection can be based on the
BIC defined by $\BIC = -2 {\cal L}_c + m \log(n)$,
where $n$ is the length of the time series and $m$ is the number of parameters estimated.
From eqn. (\theConcentratedLoglikelihood),
the loglikelihood of the $\ARzeta(i_1,\ldots,i_m)$ may be approximated by,
${\cal L}_c \approx -(n/2) \log(\hat \sigma_a^2)$.
Since
$\hat \sigma_a^2 \approx c_0 (1-\hat{\phi}_{i_1,i_1}^2) \ldots (1-\hat{\phi}_{i_m,i_m}^2)$,
where $c_0$ is the sample variance.
Hence, we obtain the approximation,
\begin{equation}
\BICzeta = \BIC(i_1,\ldots,i_m) =
n \log \left(\prod\limits_{k \in \{i_1,\ldots,i_m\}} (1-\hat{\phi}_{k,k}^2) \right) +  m \log(n).
\newcounter{BICDefinition}
\setcounter{BICDefinition}{\value{equation}}
\end{equation}
The following algorithm can be used to find the minimum $\BICzeta$ model:

\begin{description}
\item[{\it Step 1:}] Select $L$ the maximum order for the autoregression.
Select $M, M\le L$, the maximum number of parameters allowed.
The partial autocorrelation plot can be used to select $L$ large enough
so that all partial autocorrelations larger than $L$ are assumed zero.
Also, from the partial autocorrelations we can get an approximate idea of how many
partial autocorrelation parameters might be needed.
\item[{\it Step 2:}] Sort $(|\hat{\phi}_{1,1}|, \ldots, |\hat{\phi}_{L,L}|)$ in descending order
to obtain $(|\hat{\phi}_{i_1,i_1}|, \ldots, |\hat{\phi}_{i_L,i_L}|)$.
\item[{\it Step 3:}] Compute the $\BIC(i_1,\ldots,i_m)$ for $m=1,\ldots, M$ and select the minimum $\BICzeta$ model. It is usually desirable to also consider
models which are close to the minimum since sometimes these models may perform
better for forecasting on a validation sample or perhaps give better performance
on a model diagnostic check.
So in this last step, we may select the $k$ best models.
\end{description}

This polynomial time algorithm is suitable for use with long time series and with large $L$ and $M$.
Also, other criteria such as the AIC (Akaike, 1974),  $\AIC_C$ (Hurvich and Tsai, 1989)
or that of Hannan and Quinn (1979) could also be used in this algorithm.

{\noindent \bf 2.4 Residual autocorrelation diagnostics\hfill}

Let $\zeta=(\zeta_{i_1},\ldots,\zeta_{i_m})$ denote the true parameter values in an
$\ARzeta(i_1,\ldots,i_m)$ model and
let $\dot\zeta=(\dot\zeta_{i_1},\ldots,\dot\zeta_{i_m})$ denote any value in the
admissible parameter space.
Then the residuals, $\dot a_t,\ t=p+1,\ldots,n$, corresponding to the parameter $\dot \zeta$
and data $z_1,\ldots,z_n$ from a mean-corrected time series are
defined by
$\dot{a}_t = z_t - \dot{\phi}_1 z_{t-1} - \cdots - \dot{\phi}_p z_{t-p}$
where $t=p+1,\ldots,n$, $\dot \phi_i = \phi_i(\dot \zeta)$ and $p=i_m$.
The residuals corresponding to the initial values, $t=1,\ldots,p$,
may be obtained using the
backforecasting method of Box, Jenkins and Reinsel  (1994, Ch. 5) or for
asymptotic computations they can simply be set to zero.
For lag $k \ge 0$ the residual autocorrelations are defined by
$\dot{r}_k = {\dot{c}_k / \dot{c}_0},$
where $\dot{c}_k = n^{-1} \sum_{t=k+1}^n  \dot{a}_{t-k} \dot{a}_t$ for all $k \ge 0$.
When $\dot{\zeta}=\hat{\zeta}$ the residuals and residual autocorrelations will be denoted
by $\hat{a}_t$ and $\hat{r}_k$ respectively.
For any $L>1$, let
$\dot{r} = ( \dot{r}_1,\ldots,\dot{r}_L )$
and similarly for $\hat{r}$ and $r$.

{\smcaps Theorem 3.\/} $\surd n \hat r
\stackrel{D}{\longrightarrow} {\cal N}(0,{\cal V}_r)$
where
\begin{equation}
{\cal V}_r =
I_m - {\cal X} {\cal J}_\zeta {\cal I}_\zeta^{-1} {\cal J}_\zeta^{\prime} {\cal X}^{\prime}
\newcounter{Varr}
\setcounter{Varr}{\value{equation}}
\end{equation}
where $\cal X$ is the $L \times p$ matrix with $(i,j)$-entry $\psi_{i-j}$ where
the $\psi_k$ are determined as the coefficients of $B^k$ in the expansion
$1/\phi(B)=1+\psi_1 B + \psi_2 B^2 + \cdots$
and ${\cal J}_\zeta$ and ${\cal I}_\zeta$ are as defined in \S 2.3 for the
$\ARzeta(i_1,\ldots,i_m)$ model.
This theorem is proved in Appendix B.

Since ${\cal J}_\zeta^{\prime} {\cal X}^{\prime} {\cal X} {\cal J}_\zeta \approx
{\cal J}_\zeta^{\prime} {\cal I}_\phi {\cal J}_\zeta$ for $L$ large enough
and since ${\cal I}_\zeta={\cal J}_\zeta^{\prime} {\cal I}_\phi {\cal J}_\zeta$
it follows that ${\cal V}_r$ is approximately idempotent with rank $L-m$ for
$L$ large enough.
This justifies the use of the modified portmanteau diagnostic test of Ljung and Box (1978),
$Q_L = n (n+2) \sum\limits_{k=1}^L {{\hat r}_k^2 /(n-k)}$.
Under the null hypothesis that the model is adequate, $Q_L$, is approximately
$\chi^2$-distributed on $L-m$ df.

It is also useful to plot the residual autocorrelations and show their
$(1-\alpha)$\% simultaneous confidence interval.
As pointed out by Hosking and Ravishanker (1993), a simultaneous
confidence interval may be obtained by applying the Bonferonni inequality.
The estimated standard deviation of $\hat{r}_k$ is
$\EstSd(\hat{r}_k)= v_{i,i}/\surd{n}$,
where $v_{i,i}$ is the $(i,i)$ element of the covariance matrix $\hat{{\cal V}}_r$
obtained by replacing parameter $\zeta$ in eq. (\theVarr) by its estimate $\hat{\zeta}$.
Then, using the approximation with the Bonferonni inequality, it may be shown
that a $(1-\alpha)$\% simultaneous confidence interval for $\hat{r}_1,\ldots,\hat{r}_L$
is given by $\Phi^{-1}(1-\alpha/(2 m)) \EstSd(\hat{r}_k)$, where $\Phi^{-1}(\bullet)$ denotes
the inverse cumulative distribution function of the standard normal distribution.
This diagnostic plot is illustrated in \S 3.1.

\bigskip
\begin{center}
3. ILLUSTRATIVE EXAMPLES
\end{center}

{\noindent \it 3.1 Chemical process time series\hfill}

Cleveland (1971) identified an $\ARphi(1,2,7)$
and Unnikrishnan (2004) identified an $\ARphi(1,3,7)$ model
for Series A (Box, Jenkins and Reinsel, 1994).
Either directly from the partial autocorrelation plot in Figure 1
or using the $\BICzeta$ algorithm in \S 2.3 with $L=20$ and $M=10$,
an $\ARzeta(1,2,7)$ subset is selected.
The top five models with this algorithm are shown in Table I.
Figure 2 shows the residual autocorrelation plots for the fitted $\ARphi$ and $\ARzeta$ models.
The respective maximized loglikelihoods were ${\cal L}_c = 232.96$ and $229.42$ respectively.
Thus, a slightly better fit is achieved by the $\ARphi$ model in this case, but
since the difference is small, it may be concluded that both models fit about equally well.

{\noindent \it 3.2 Forecasting experiment\hfill}

McLeod and Hipel (1995) fitted an $\ARphi(1,9)$ to the treering series identified
as {\smcaps Ninemile} in their article.
This series consists of $n=771$ consecutive annual treering width measurements
on Douglas fir at Nine Mile Canyon, Utah for the years 1194$-$1964.
For our forecasting experiment the first 671 values were used as the training data
and the last 100 as the test data.
The partial autocorrelation plot of the training series is shown in Figure 3.
This plot suggests $L=20$ and $M=10$ in the algorithm in \S 2.3 will suffice.
The three best $BICzeta$ models were $\ARzeta(1)$, $\ARzeta(1,9)$ and $\ARzeta(1,2,9)$.
After fitting with exact maximum likelihood,
the one-step forecast errors were computed for the test data.
The $\ARphi(1,9)$ model was also fit to the training series and
the one-step forecast errors over the next 100 values were computed.
Table II compares the fits achieved as well as the root mean square error on the test data.
From Table II as well as with further statistical tests, it was concluded that there is
no significant difference in forecast performance.
\medskip

{\noindent \it 3.3 Monthly sunspot series\hfill}

The monthly sunspot numbers, 1749$-$1983 (Andrews and Hertzberg, 1985),
are comprised of $n=2820$ consecutive values.
Computing the coefficient of skewness for the transformed data, $z_t^\lambda$,
with $\lambda = 1, 2/3, 1/2, 1/3$ we obtained
$g_1 =  1.10, 0.48, 0.09,-0.45$ respectively.
It is seen that a square-root transformation will improve the normality assumption.
For the square-root transformed series, subset $\ARzeta$ models were determined
using the $\AICzeta$ and $\BICzeta$ algorithms with $L=300$
and $M=100$.
These algorithms produced subset models with $m=70$ and $m=20$ autoregressive
coefficients.
Maximum likelihood estimation of these two models only required about
30 minutes and 3 minutes respectively on a 3 GHz PC using our {\it Mathematica\/} software.
The best nonsubset $\AR(p)$ models for the square-root monthly sunspots
using the AIC and BIC are compared with the subset models in Table III.
The $\ARzeta$ fitted with the $\BICzeta$ algorithm has fewer parameters than each of these nonsubset models and it performs better on both the AIC and BIC criteria.
Residual autocorrelation diagnostic checks did not indicate any model inadequacy in any of the
fitted models.

\bigskip
\begin{center}
4. CONCLUDING REMARKS
\end{center}
\medskip

The methods presented in this paper can be extended to subset MA models.
In this case the $\zeta$-parameters are inverse partial autocorrelations
(Hipel and McLeod, 1994, \S 5.3.7).
Bhansali (1983) showed that the distribution of the inverse partial
autocorrelations is equivalent to that of the partial autocorrelations,
so the model selection using a modified inverse partial autocorrelation plot
or AIC/BIC criterion may be implemented.
Similarly the distribution of the residual autocorrelations is essentially
equivalent to that given in our Theorem 3.

As discussed by Monahan (1984) and Marriott and Smith (1992) the transformation
used in eq. (\theBTransformation) may be extended to reparameterize ARMA models.
Hence the subset AR model may be generalized in this way to the subset ARMA case.

Software written in {\it Mathematica\/} is available from the authors'
web page for reproducing the examples given in this article
as well as for more general usage.

\newpage
\bigskip
\begin{center}
APPENDIX A: PROOF OF THEOREM 1
\end{center}
\medskip

Let $\dot{\zeta}$, $\hat{\zeta}$ and $\zeta$ denote respectively a vector of parameters
in the admissible region, the maximum likelihood estimate and the true parameter
values and
similarly for other functions of these quantities such as the likelihood and residuals.
Without loss of generality we may assume that $\mu=0$ and $\sigma=1$ are known.
Ignoring terms which are $O_p(1)$, the loglikelihood function of $\dot{\zeta}$
may be written,
\begin{equation}
\dot{\cal L}=-{1\over 2} \sum \dot{a}_t^2,
\newcounter{LoglikelihoodAppendixA}
\setcounter{LoglikelihoodAppendixA}{\value{equation}}
\end{equation}
where $\dot{a}_t = z_t - \dot{\phi}_1 z_{t-1} - ... - \dot{\phi}_p z_{t-p}$.
Note that for all $i$ and $j$,
\begin{equation}
{\partial \dot{a}_t \over \partial \dot{\phi}_i}=-z_{t-i}
\newcounter{FirstDerivativea}
\setcounter{FirstDerivativea}{\value{equation}}
\end{equation}
and
\begin{equation}
{\partial \dot{a}_t \over {\partial \dot{\phi}_i \dot{\phi}_j}}=0.
\newcounter{SecondDerivativea}
\setcounter{SecondDerivativea}{\value{equation}}
\end{equation}
It follows that
\begin{equation}
{1\over n}{\partial {\cal L} \over \partial \zeta} =
{1\over n} \sum a_t {\cal J}_{\zeta}^{\prime} (z_{t-1},...,z_{t-p})^{\prime},
\newcounter{LoglikelihoodDerivative}
\setcounter{LoglikelihoodDerivative}{\value{equation}}
\end{equation}
where $(z_{t-1},...,z_{t-p})^{\prime}$ denotes the transpose of the $p$-dimensional
row vector and
\begin{equation}
{\cal J}_\zeta ={\partial (\phi_1,\ldots,\phi_p) \over \partial (\zeta_{i_1},\dots,\zeta_{i_m})}.
\newcounter{JacobianZetaAppendixA}
\setcounter{JacobianZetaAppendixA}{\value{equation}}
\end{equation}
Since $\sum z_{t-j} a_t /n \stackrel{p} {\longrightarrow} 0$ when $j>0$,
it follows that from (\theLoglikelihoodDerivative),
\begin{equation}
{n^{-1}}{\partial {\cal L} \over \partial \zeta} \stackrel{p} {\longrightarrow} 0.
\newcounter{ScoreFunctionAppendixA}
\setcounter{ScoreFunctionAppendixA}{\value{equation}}
\end{equation}
Similarly, neglecting terms which are $O_p(1/{\surd n})$,
\begin{eqnarray}
-{1\over n} {\partial^2 {\cal L} \over \partial \zeta \partial \zeta^{\prime}}
& = &
{1\over n} {\cal J}_\zeta^{\prime} (\sum z_{t-i} z_{t-j}) {\cal J}_\zeta
\nonumber\\
& \stackrel{p} {\longrightarrow} &
{\cal J}_\zeta^{\prime} {\cal I}_\phi {\cal J}_\zeta
\nonumber\\
& = &
{\cal I}_{\zeta},
\newcounter{InformationMatrixZetaAppendixA}
\setcounter{InformationMatrixZetaAppendixA}{\value{equation}}
\end{eqnarray}
where
\begin{equation}
{\cal I}_\phi=\left(\Cov(z_{t-i}, z_{t-j})\right)_{p \times p}.
\newcounter{InformationMatrixARAppendixA}
\setcounter{InformationMatrixARAppendixA}{\value{equation}}
\end{equation}
Since the transformation is one-to-one, continuous and differentiable,
it follows that ${\cal I}_\zeta$ must be positive definite since
${\cal I}_\phi$ is positive definite.
Expanding
${n^{-1} \partial \dot{\cal L} / \partial \dot{\zeta}}$
about $\dot{\zeta}=\hat{\zeta}$ and evaluating at
$\dot{\zeta}=\zeta$ and noting that third and higher-order terms are zero,
\begin{equation}
0={1\over n} {\partial {\cal L} \over \partial \zeta}
+
(\hat{\zeta}-\zeta){1\over n} {\partial^2 {\cal L} \over \partial \zeta \partial \zeta^{\prime}}.
\newcounter{LoglikelihoodExpansion}
\setcounter{LoglikelihoodExpansion}{\value{equation}}
\end{equation}
It follows from eq. (\theScoreFunctionAppendixA) and (\theLoglikelihoodExpansion) that
$\hat \zeta \stackrel{p} {\longrightarrow} \zeta$.
Since
\begin{equation}
{1\over n} {\partial^2 {\cal L} \over \partial \zeta \partial \zeta^{\prime}}
=
{\cal I}_{\zeta}+O_p({1 \over \surd{n}}),
\newcounter{LoglikelihoodExpansionPartTwo}
\setcounter{LoglikelihoodExpansionPartTwo}{\value{equation}}
\end{equation}
it follows that
\begin{equation}
-{1\over \surd{n}} {\partial {\cal L} \over \partial \zeta}
=
\surd{n} ( \hat{\zeta}-\zeta ) {\cal I}_{\zeta}+o_p(1).
\newcounter{LoglikelihoodExpansionPartThree}
\setcounter{LoglikelihoodExpansionPartThree}{\value{equation}}
\end{equation}
Since
\begin{equation}
{1\over \surd{n}} {\partial {\cal L} \over \partial \phi}
\stackrel{D} {\longrightarrow}
N(0, {\cal I}_\phi).
\newcounter{ScoreConvergenceAR}
\setcounter{ScoreConvergenceAR}{\value{equation}}
\end{equation}
and
\begin{equation}
{\partial {\cal L} \over \partial \zeta}
=
{\partial {\cal L} \over \partial \phi}^{\prime} {\cal J}_{\zeta}
\end{equation}
it follows that
\begin{equation}
\surd{n} ( \hat{\zeta}-\zeta )
\stackrel{D} {\longrightarrow}
N(0, {\cal I}_\zeta^{-1}).
\newcounter{SubsetDistribution}
\setcounter{SubsetDistribution}{\value{equation}}
\end{equation}

\newpage
\bigskip
\begin{center}
 APPENDIX B: PROOF OF THEOREM 3
\end{center}
\medskip

Without loss of generality and ignoring terms which are $O_p(1)$
we may write the loglikelihood function as
${\cal L}(\dot \zeta) = -{1 \over 2} \sum \dot{a}_t^2$.
Then $r_k = c_k + O_p(1/n)$
and,
\begin{equation}
\zeta - \hat{\zeta} = {\cal I}_\zeta^{-1} s_c + O_p(1/n),
\newcounter{ZetaEfficientScore}
\setcounter{ZetaEfficientScore}{\value{equation}}
\end{equation}
where
\begin{eqnarray}
s_c & = & {\partial {\cal L} \over \partial \zeta} \nonumber\\
& = & {\cal J}_\zeta {\partial {\cal L} \over \partial \phi} \nonumber\\
& = & {\cal J}_\zeta \left( \sum a_t z_{t-i} \right).
\newcounter{EfficientScore}
\setcounter{EfficientScore}{\value{equation}}
\end{eqnarray}
It follows that, neglecting a term which is $O(1)$,
\begin{equation}
n \Cov(\hat{\zeta}, r)=-{\cal I}_\zeta^{-1} {\cal J}_\zeta^{\prime} {\cal X}^{\prime}.
\newcounter{CrossCovariance}
\setcounter{CrossCovariance}{\value{equation}}
\end{equation}
Expanding $\dot{r}$ about $\dot{\zeta}=\zeta$ and evaluating at
$\dot{\zeta}=\hat{\zeta}$,
\begin{equation}
\hat{r} = r + {\cal X} {\cal J}_\zeta (\hat{\zeta}-\zeta).
\newcounter{rExpansion}
\setcounter{rExpansion}{\value{equation}}
\end{equation}
From (\therExpansion), it follows
that $\surd{n} r$ is asymptotically normal with mean zero and
the given covariance matrix.

This theorem could also be derived using the result of Ahn (1988) on
multivariate autoregressions with structured parameterizations.

\newpage
\begin{center}
{\smcaps REFERENCES\hfill}
\end{center}
\parindent 0pt

\hind {\smcaps Ahn, S.K.\/} (1988)
Distribution for residual autocovariate in multivariate autoregressive models
with structured parameterization.
{\it Biometrika\/} 75, 590--593.

\hind
{\smcaps Akaike, H.\/} (1974)
A new look at the statistical model identification.
{\it IEEE Transactions on Automatic Control\/} AC-19, 716--723.

\hind {\smcaps Andrews, D.F. and Herzberg, A.M.\/} (1985)
{\it Data: A Collection of Problems from Many Fields for the Student and Research Worker\/}.
New York: Springer-Verlag.

\hind
{\smcaps Barndorff-Nielsen, O. and Schou G.\/} (1973)
On the parametrization of autoregressive models by partial autocorrelations.
{\it Journal of Multivariate Analysis\/}  3, 408--419.

\hind
{\smcaps Bhansali, R.J.\/} (1983)
The inverse partial autocorrelation function of a time series and its applications.
{\it Journal of Multivariate Analysis\/} 13, 310--327.

\hind
{\smcaps Box, G.E.P., Jenkins, G.M. and Reinsel, G.C.\/}  (1994)
{\it Time Series Analysis: Forecasting and Control\/},
3rd Ed., San Francisco: Holden-Day.

\hind
{\smcaps Box, G.E.P. and Luce\~no, A.\/} (1997)
{\it Statistical Control by Monitoring and Feedback Adjustment\/},
New York: Wiley.

\hind
{\smcaps Chen, C.W.S.\/} (1999)
Subset selection of autoregressive time series models.
{\it Journal of Forecasting\/} 18, 505--516.

\hind
{\smcaps Cleveland, W.S.\/} (1971)
The inverse autocorrelations of a time series and their applications.
{\it Technometrics\/} 14, 277--298.

\hind
{\smcaps George, E.I. and McCulloch, R.E.\/} (1993)
Variable selection via Gibbs sampling.
{\it Journal of the American Statistical Association\/} 88, 881-889.

\hind
{\smcaps Haggan, V. and Oyetunji, O.B.\/}  (1984)
On the selection of subset autoregressive time series models.
{\it Journal of Time Series Analysis\/} 5, 103--113.

\hind
{\smcaps Hannan, E.J. and Quinn, B.G.\/} (1979)
The determination of the order of an autoregression.
{\it Journal of the Royal Statistical Society\/} B 41, 190--195.

\hind
{\smcaps Hipel, K.W. and McLeod, A.I.\/} (1994)
{\it Time Series Modelling of Water Resources and Environmental Systems,\/}
Amsterdam: Elsevier.

\hind
{\smcaps Hosking, J.R.M. and Ravishanker, N.\/} (1993)
Approximate simultaneous significance intervals for residual autocorrelations of
autoregressive-moving average time series models.
{\it Journal of Time Series Analysis\/} 14, 19--26.

\hind
{\smcaps Hurvich, M.C. and Tsai, C.\/} (1989)
Regression and time series model selection in small samples.
{\it Biometrika\/} 76, 297--307.

\hind
{\smcaps Ljung, G.M. and Box, G.E.P.\/} (1978)
On a measure of lack of fit in time series models.
{\it Biometrika\/} 65, 297--303.

\hind
{\smcaps Marriott, J.M. and Smith, A.F.M.\/} (1992)
Reparameterization aspects of numerical Bayesian methodology for
autoregressive-moving average models.
{\it Journal of Time Series Analysis\/} 13, 327--343.

\hind
{\smcaps McClave, J.\/} (1975)  Subset autoregression.
{\it Technometrics\/} 17, 213--220.

\hind
{\smcaps McLeod, A.I. and Hipel, K.W.\/} (1995)
Exploratory spectral analysis of hydrological time series.
{\it Journal of Stochastic Hydrology and Hydraulics\/} 9, 171-205.

\hind
{\smcaps McLeod, A.I.\/} (1978)
On the distribution of residual autocorrelations in Box-Jenkins models.
{\it Journal of the Royal Statistical Society B\/} 40, 396--402.

\hind
{\smcaps Monahan, J.F.\/} (1984)
A note on enforcing stationarity in autoregressive-moving average models.
{\it Biometrika\/} 71, 403--404.

\hind
{\smcaps Percival, D.B. and Walden, A.T.\/} (1993)
{\it Spectral Analysis For Physical Applications\/},
Cambridge, Cambridge University Press.

\hind
{\smcaps Siddiqui, M.M.\/} (1958)
On the inversion of the sample covariance matrix in a stationary autoregressive process.
{\it Annals of Mathematical Statistics\/} 29, 585--588.

\hind
{\smcaps Taniguchi, M.\/} (1983)
On the second order asymptotic efficiency of estimators of gaussian ARMA processes.
{\it The Annals of Statistics\/} 11, 157--169.

\hind
Tj\o stheim, D. \& Paulsen, J. (1983)
Bias of some commonly-used time series estimates.
{\it Biometrika\/} {\bf 70,\/} 389--399.
Correction {\it Biometrika\/} {\bf 71,\/} p. 656.

\hind
{\smcaps Tong, H.\/} (1977)
Some comments on the Canadian lynx data.
{\it Journal of the Royal Statistical Society A\/} 140, 432--436.

\hind
{\smcaps Unnikrishnan, N. K.\/} (2004)
Bayesian Subset Model Selection for Time Series.
{\it Journal of Time Series Analysis\/} 25, 671-690

\hind
{\smcaps Zhang, Y. and McLeod, A.I.\/} (2005, under review)
Computer algebra derivation of the bias of Burg estimators.

\hind
{\smcaps Zhang, X. and Terrell, R.D.\/} (1997)
Prediction modulus: A new direction for selecting subset autoregressive models.
{\it Journal of Time Series Analysis\/} 18, 195--212.

\newpage
\begin{table}
{\smcaps Table I.
Top five models for Series A using the $\BICzeta$ algorithm.
\/
}
$$
\vbox{\halign{\strut
#\quad\hfill    
&\hfill#\quad\hfill    
\cr             
\noalign{\hrule}
\noalign{\smallskip}
\noalign{\hrule}
\noalign{\smallskip}
Model&$\BICzeta$\cr
\noalign{\smallskip}
\noalign{\hrule}
\noalign{\smallskip}
$\ARzeta(1,2,7)$&$-82.2$\cr
$\ARzeta(1,2,7,15)$&$-81.5$\cr
$\ARzeta(1,2)$&$-80.4$\cr
$\ARzeta(1,2,7,6,15)$&$-80.4$\cr
$\ARzeta(1,2,7,6,15,17)$&$-78.2$\cr
\noalign{\smallskip}
\noalign{\hrule}
}}$$
\end{table}
\strut
\vfill\eject\newpage

\newpage
\begin{table}
{\smcaps Table II.
Comparison of models fit to the training portion of the Ninemile dataset
and the root mean square error, RMSE, of the one step-ahead forecasts on the
test portion.\/
}
$$\vbox{\halign{\strut
#\quad\hfill    
&\hfill#\quad\hfill    
&\hfill#\quad\hfill    
&\hfill#\quad\hfill    
&\hfill#\quad\hfill    
\cr             
\noalign{\hrule}
\noalign{\smallskip}
\noalign{\hrule}
\noalign{\smallskip}
Model&${\cal L}_c$&BIC&$\hat \sigma_a$&RMSE
\cr
\noalign{\smallskip}
\noalign{\hrule}
\noalign{\smallskip}
$\ARphi(1,9)$&$-2462.9$&$4938.8$&$39.3$&$42.3$
\cr
$\ARzeta(1)$&$-2467.9$&$4942.4$&$39.6$&$43.3$
\cr
$\ARzeta(1,9)$&$-2465.3$&$4943.6$&$39.4$&$42.8$
\cr
$\ARzeta(1,2,9)$&$-2463.0$&$4945.4$&$39.3$&$42.3$
\cr
\noalign{\smallskip}
\noalign{\hrule}
}}$$
\end{table}
\strut
\vfill\eject\newpage

\begin{table}
{\smcaps Table III.
Comparison of models fitted to the square-root of the monthly sunspot series.
The $\ARzeta$ model fitted using the AIC is denoted by $\ARzeta$(AIC) and
similarly for $\ARzeta$(BIC).
The best nonsubset AR models fitted using the AIC and BIC are
denoted by $\AR$(AIC) and $\AR$(BIC).
The number of coefficients, $m$, is shown.
The values of the AIC and BIC shown are obtained using the exact loglikelihood.
}
$$\vbox{\halign{\strut
#\quad\hfill    
&\hfill#\quad\hfill    
&\hfill#\quad\hfill    
&\hfill#\quad\hfill    
&\hfill#\quad\hfill    
\cr             
\noalign{\hrule}
\noalign{\smallskip}
\noalign{\hrule}
\noalign{\smallskip}
Model&$m$&${\cal L}_c$&AIC&BIC
\cr
\noalign{\smallskip}
\noalign{\hrule}
\noalign{\smallskip}
$\AR_\zeta$(AIC)&$70$&$-148.2$&$436.4$&$852.6$
\cr
$\AR_\zeta$(BIC)&$20$&$-236.5$&$513.$&$631.9$
\cr
$\AR$(AIC)&$27$&$-241.1$&$536.3$&$696.8$
\cr
$\AR$(BIC)&$21$&$-252.5$&$547.$&$671.8$
\cr
\noalign{\smallskip}
\noalign{\hrule}
}}$$
\end{table}

\newpage
\strut
\begin{figure}
\centerline{\epsfig{figure=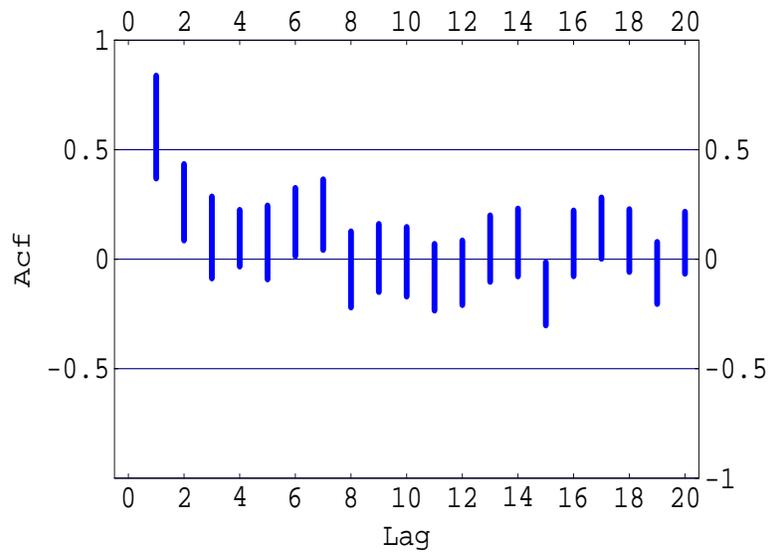,width=4in,height=3in}}
\caption{Partial autocorrelation plot of Series A.}
\end{figure}

\newpage
\strut
\begin{figure}
\centerline{\epsfig{figure=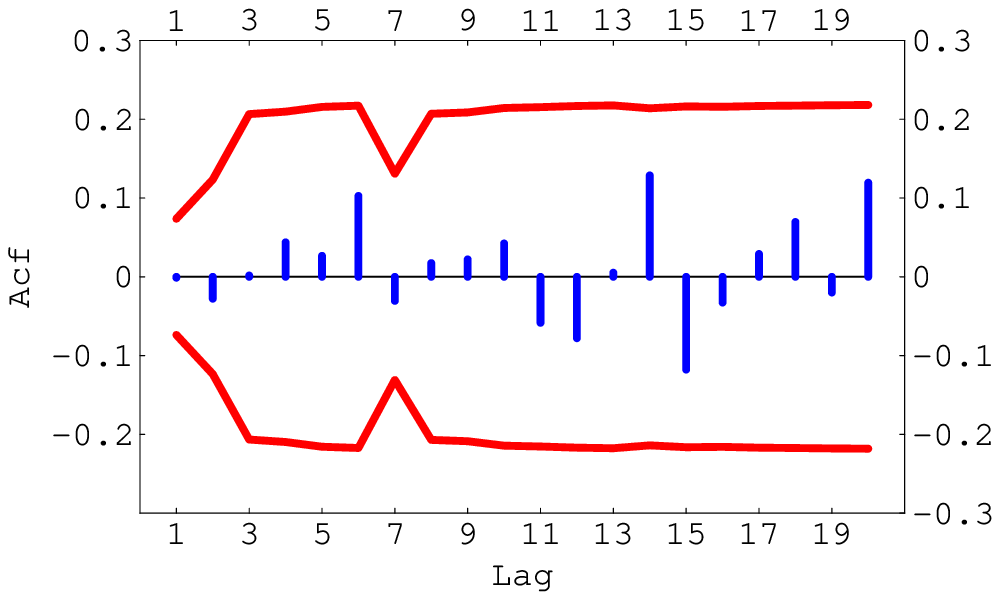,width=4in,height=3in}}
\centerline{\epsfig{figure=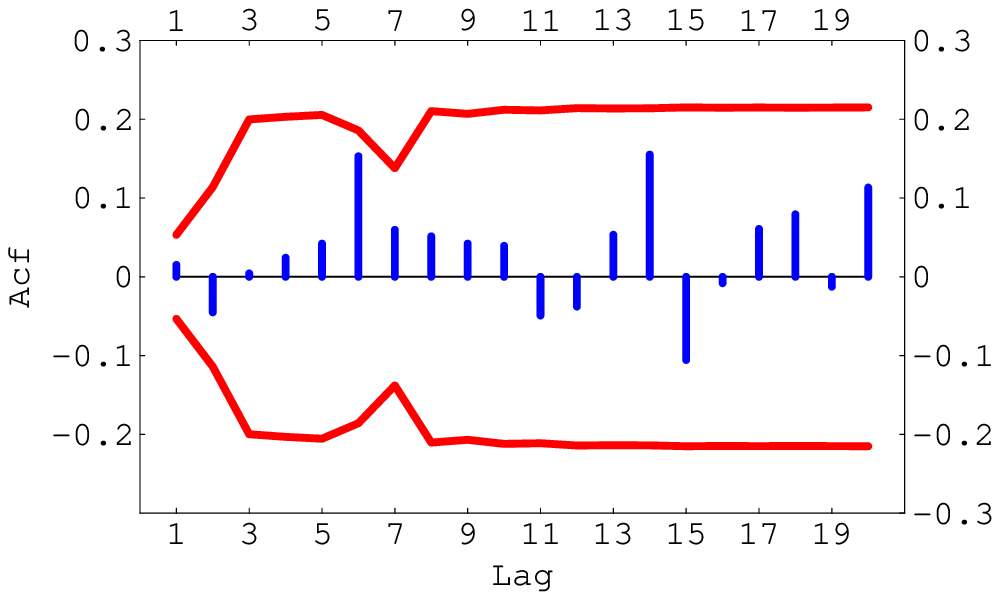,width=4in,height=3in}}
\caption{Residual autocorrelation plots for the $\ARphi(1,2,7)$, upper panel, and
the $\ARzeta(1,2,7)$, lower panel.}
\end{figure}

\newpage
\strut
\begin{figure}
\centerline{\epsfig{figure=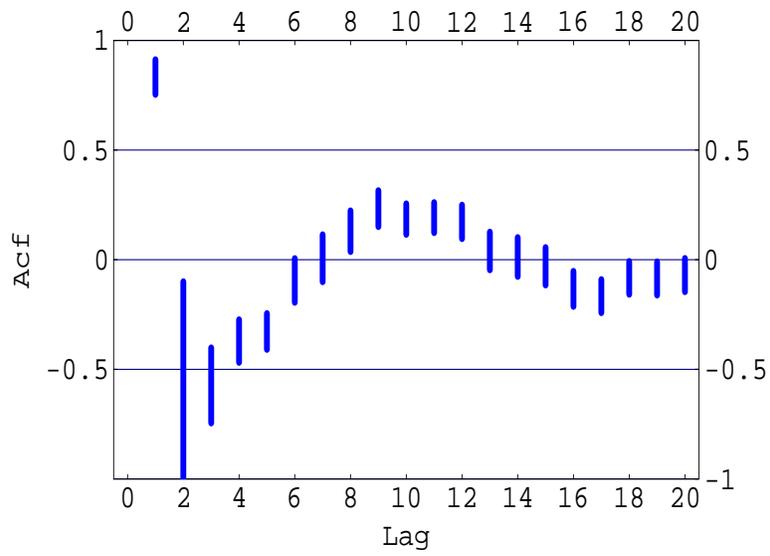,width=4in,height=3in}}
\caption{Partial autocorrelation plot of the training portion of the Ninemile
treering series.}
\end{figure}

\end{document}